\documentclass[12pt]{article}

\usepackage{flushend}
\usepackage{balance}
\usepackage[utf8]{inputenc}
\usepackage[english]{babel}
\usepackage{amsmath}
\usepackage{amsfonts}
\usepackage{amssymb}
\usepackage{epsfig}
\usepackage{enumerate}
\usepackage{fleqn}
\usepackage{color}
\usepackage{url}
\usepackage{mathtools}
\usepackage{graphicx}
\usepackage{hyperref}
\usepackage{subfigure}
\usepackage[draft]{changes} 
\usepackage{lineno}  

\newtheorem{assumption}{Assumption}

\newtheorem{theorem}{Theorem}

\newtheorem{proposition}{Proposition}

\usepackage{times}

\newenvironment{sciabstract}{%
\begin{quote} \bf}
{\end{quote}}

\title{Lightning optimizes: a threshold mechanism\\ ensures minimum-path flow}

\author
{Franco Blanchini$^{1\ast}$, Daniele Casagrande$^{2}$,  Filippo Fabiani$^{3}$,\\
Giulia Giordano$^{4}$,  David Palma$^{1}$, Raffaele Pesenti$^{5}$
\\
 \normalsize{$^{1}$Dipartimento di Matematica, Informatica e Fisica, Universit\`a di Udine, 33100 Udine, Italy }\\
\normalsize{$^{2}$Dipartimento Politecnico di Ingegneria e Architettura, Universit\`a di Udine, 33100 Udine, Italy }\\
\normalsize{$^{3}$Department of Engineering Science, University of Oxford, OX1 3PJ, United Kingdom }\\
\normalsize{$^{4}$Dipartimento di Ingegneria Industriale, Universit\`a di Trento, 38123 Povo (TN), Italy}\\
\normalsize{$^{5}$Dipartimento di Management, Universit\`a Ca' Foscari, 30121 Venezia, Italy}\\
\normalsize{$^\ast$To whom correspondence should be addressed; E-mail:  blanchini@uniud.it.}
}

\date{}

\begin{document} 


\baselineskip24pt


\maketitle 


\begin{sciabstract}
  A well-known property of linear resistive electrical networks is that the current distribution minimizes the total dissipated energy.
When the circuit includes resistors with nonlinear monotonic characteristic, the current distribution minimizes in general a different functional.
We show that, if the nonlinear characteristic is a threshold-like function and the energy generator is concentrated in a single point, as in the case of lightning or dielectric discharge, then the current flow is concentrated along a single path, which is a minimum path to the ground with respect to the threshold.
We also propose a dynamic model that explains and qualitatively reproduces the lightning transient behavior: initial generation of several plasma branches and subsequent dismissal of all branches but the one reaching the ground first, which is the optimal one.
\end{sciabstract}
 
\section*{Introduction}

In lightning or gas electrical discharge, the current flow is essentially concentrated along a single path. Under very slow motion it can be seen that lightning starts by generating several branches and then develops by dismissing all of them but a single one, along which the energy is discharged \cite{video}.
This phenomenon has been deeply investigated. Several types of lightning are known, which are carefully described, e.g., in \cite{Cooray2008,Uman2001}. 
As far as the numerical modeling of the phenomenon is concerned, computational models for lightning simulation have been proposed in \cite{DeConti2008,Hager1998,Podgorski1987,Sarajcev2008}, while the fractal nature of  lightning discharge has been investigated in \cite{Niemeyer1984,Nguyen2001,Sanudo1995,Theethayi2005}.
Detailed surveys on the subject are also available; see, e.g.,~ \cite{Rakov2009} and \cite{Rakov1998}.

Here, we do not investigate the whole phenomenon in its complexity. We rather focus on a specific question about path formation in lightning discharge: we are interested in the initial phase of the process, when the lightning path is formed.
Also, we consider the ideal case in which the lightning source is a single point and the
final destination is a zero-potential ground. This type of lightning, classified as Category 1 Lightning, ``is the most common cloud-to-ground lightning. It accounts for over 90\% of the worldwide cloud-to-ground flashes'' \cite{Uman2001}.
Cloud-to-ground lightning begins with an initial breakdown and the consequent creation of a ionized channel, the \emph{stepped leader}, which generates several branches. Once the stepped leader is close to the ground, it may be approached by channels originating from the ground, the \emph{connecting leaders}. When the stepped leader finally connects the ground to the cloud, the \emph{return stroke} is triggered, which is a ground-potential upward wave \cite{Cooray2008,Visacro}.
After the return stroke reaches back the cloud, the main branch reaching
the ground is crossed by a long-duration discharging current: the \emph{continuing currents}. In the meanwhile, the secondary branches originally established
by the stepped leader are depleted: the continuing currents flow along
the main path only \cite{Cooray2008}.

Why is the lightning current eventually concentrated along a single path? Does this path enjoy any optimality property?

To mathematically address these questions, we consider an idealized model based on the assumption that lightning is mainly due to a dielectric breakdown of the air (gas in the case of discharges). The current-voltage diagram of a gas is characterized by two regions: for all the voltage values belonging to a symmetric interval around the origin, the current is very low (high resistivity); for voltage values outside this interval, the current becomes very large (low resistivity). The voltage value corresponding to the ends of the interval is called \emph{breakdown threshold}.
Then, we consider an ideal characteristic with conductivity approaching infinity when 
the electric field is larger than a threshold \cite{Hager1989,Hager1998}.

Lightning path can be interpreted as the solution of an optimization problem over a network.
To formulate the problem, we consider a graph describing an electrical network, where capacitances and resistors with possibly nonlinear characteristic are associated with the links. The grid model we use is akin to that proposed in \cite[Eq. 5]{Hager1998}, which is the discretized version of nonlinear field equations \cite{Hager1989,Hager1998}.
We show that the steady-state solution minimizes a convex functional that, in the special case of linear resistors, turns out to be the total dissipated energy.
Conversely, if the resistor characteristic is a threshold-like function, the steady-state
solution becomes the \emph{minimum path}, where each current link is weighted by its local threshold voltage; hence, the optimal path is the ``minimum-threshold path''. 
Our main result is supported by a theorem and reinforced by simulations of randomly generated graphs with random threshold values, showing that the transient behavior of the model can faithfully reproduce the qualitative lightning evolution (see Materials and Methods).


\begin{figure} [h]
		\centering
		\includegraphics[width=0.5\columnwidth]{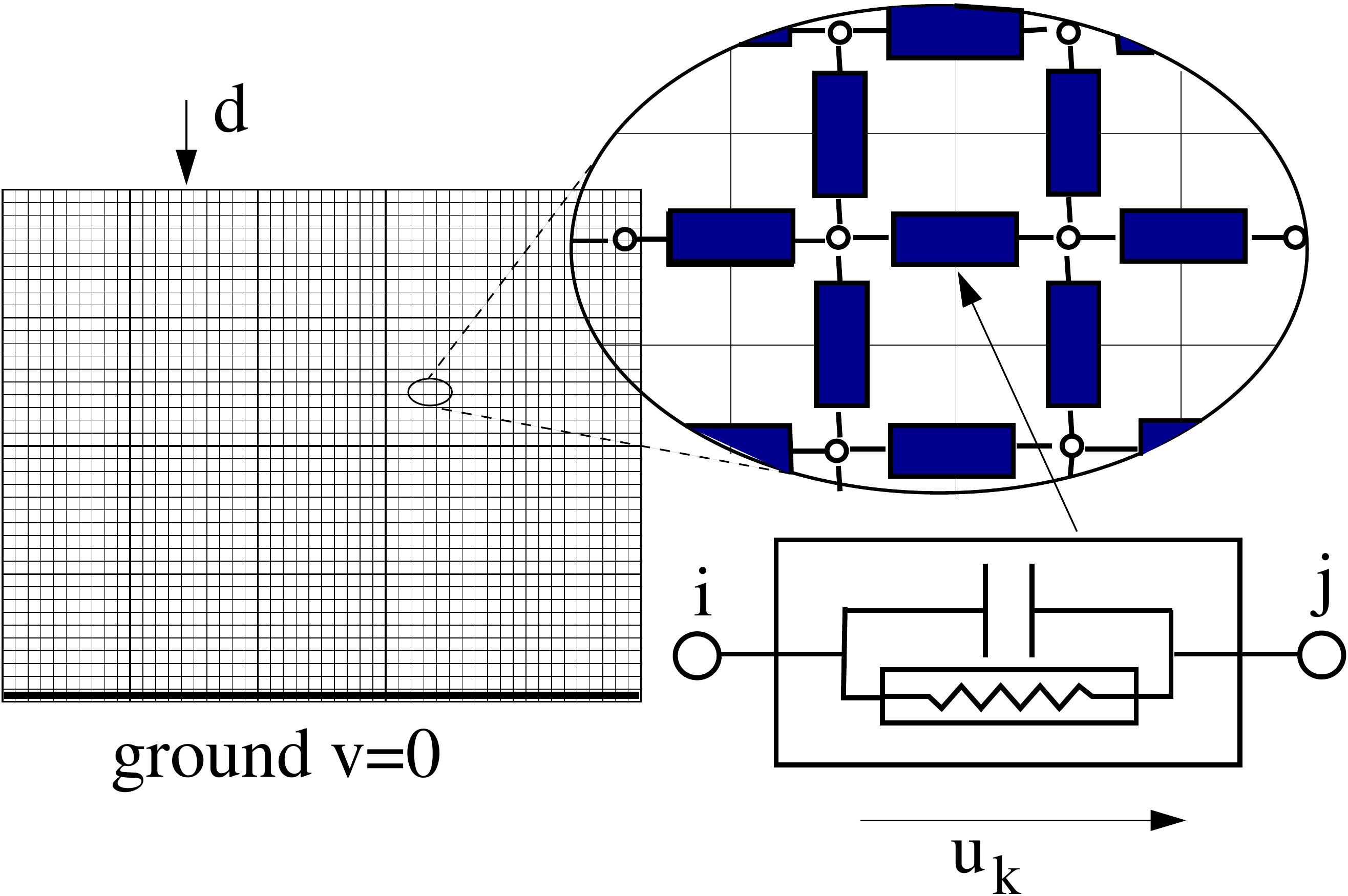}
		\caption{The electrical grid model with capacitances and possibly nonlinear resistors connecting adjacent terminals. In the graph representation, each terminal is associated with a node and each electric component with a link.}
		\label{fig:circuit}
\end{figure}

\section*{Results}

So, how does lightning optimize its path?
To build a model, we consider the electrical grid network in Fig.~\ref{fig:circuit}, with a capacitive and a possibly nonlinear resistive effect between adjacent terminals.
Ground terminals are connected among them with zero resistance (ideal conductor ground) and the ground voltage is $v=0$.
The network is associated with a graph where the $n+1$ nodes represent the terminals and the links represent the electric impedances. 
In particular, node $n$ corresponds to the zero-potential ground, while at the source node $0$ a (current or voltage) generator is applied, with its other terminal grounded, inducing an input current $d$ that enters the network.

\begin{figure} [h!]
		\centering
		\includegraphics[width=0.25\columnwidth]{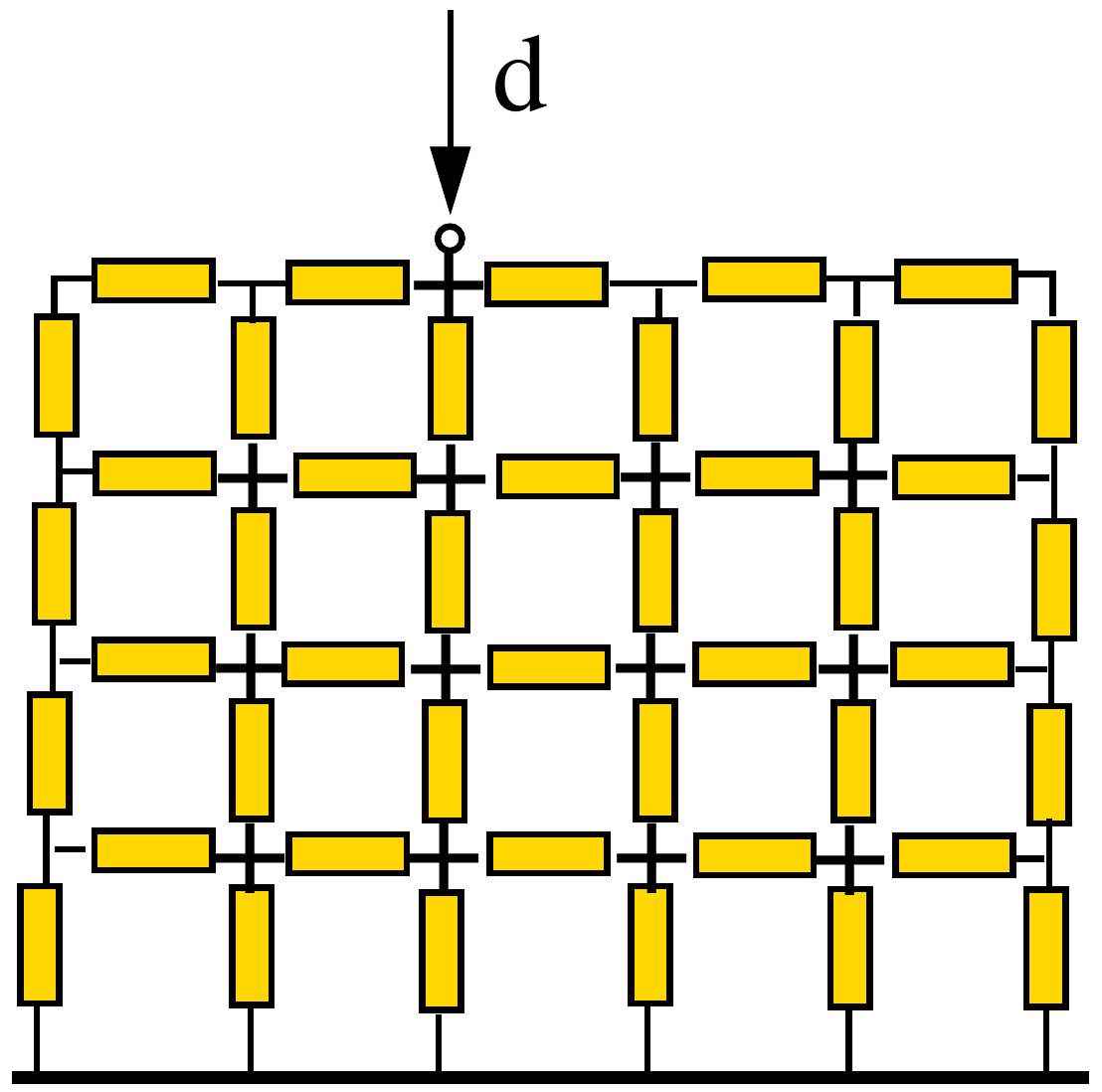}~~~~
                \includegraphics[width=0.25\columnwidth]{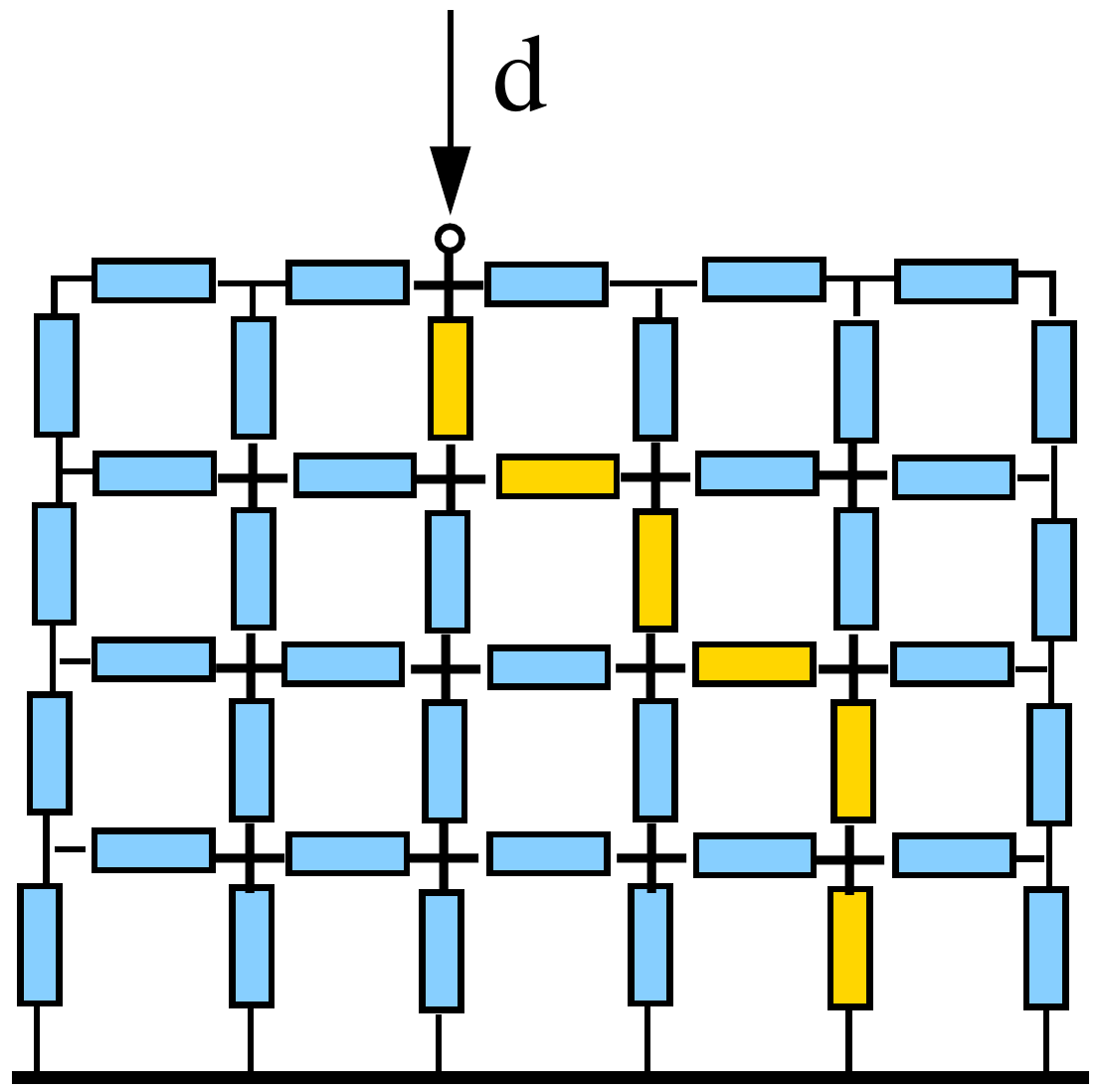}
		\caption{The electrical current distribution at steady state:  
scattered in the case of linear resistances (left); concentrated along 
a single path in the case of threshold characteristics (right), when the current flows along the path that minimizes the sum of the dielectric rigidities of its links. Yellow (resp. blue) means presence (resp. absence) of flowing current.}
		\label{fig:scat_conc}
\end{figure}

As shown in Fig.  \ref{fig:circuit}, each link is assumed to be the parallel connection of a capacitor and a 
possibly nonlinear resistor. Injecting a current $d$ in a node of the network
leads, after a transient, to a steady state in which the currents flow only through the resistors. 
If these are linear,
the steady-state solution corresponds to the current distribution that minimizes the total dissipated energy \cite{AMO93}:
\begin{equation}
E_{tot}=\sum_{k} R_k \bar u_k^2\, \rightarrow \min\,, \label{eq:dissipated}
\end{equation}
where $R_k$ is the resistance value associated with link $k$ and $\bar u_k$ is the steady-state current 
flowing through it.
In this minimum energy configuration, steady-state currents are scattered all over the network (as in Fig. \ref{fig:scat_conc}, left).
In phenomena such as lightning, the situation is completely different \cite{Cooray2008,Uman2001}:
after a transient, lightning ``chooses" a single path (as in Fig. \ref{fig:scat_conc}, right). Why and how is this single path chosen?

Assume that, denoting by $v_i$ and $v_j$ the potentials at the extreme nodes of link $k$, the resistance obeys the nonlinear law
\begin{equation*}
\begin{cases}
R_k^{th}=u_k/(v_i-v_j) = \infty, & \mbox{if} \,\, |v_i-v_j| < V_k,\\
R_k^{th}=u_k/(v_i-v_j) = 0, & \mbox{if} \,\, |v_i-v_j|\geq V_k,
\end{cases} 
\end{equation*}
where the value of the threshold $V_k$ is the \emph{local} dielectric rigidity of link $k$.

				\begin{figure} [h!]
		\centering
		\includegraphics[width=0.5\columnwidth]{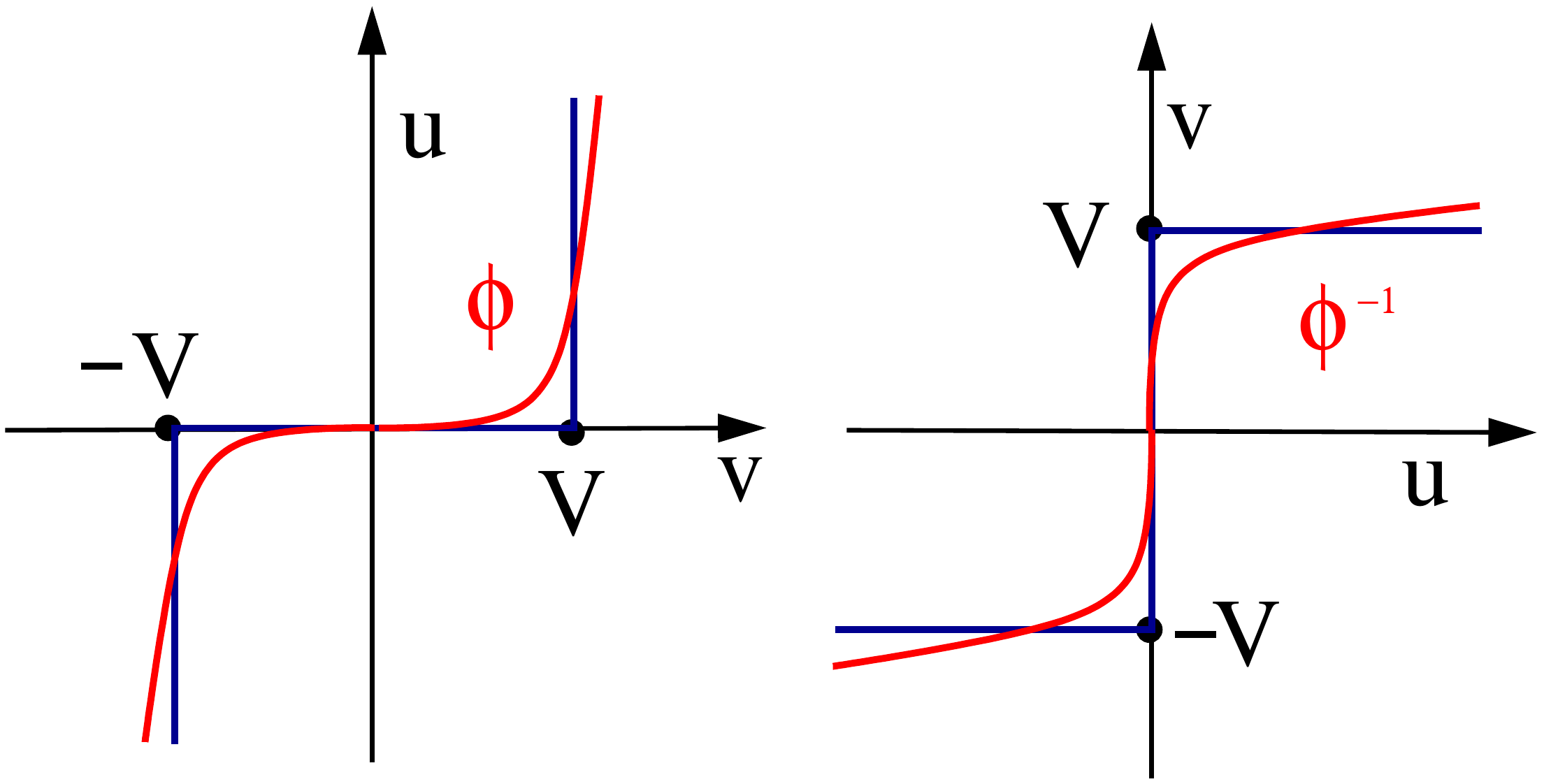}  
		\caption{Nonlinear current-voltage and voltage-current characteristics, where $v$ denotes voltage and $u$ current. A generic threshold-like current-voltage characteristic $\phi$, with threshold $V$ (red, left) and its inverse $\phi^{-1}$ (red, right). The ideal, sharp threshold characteristics (blue) can be seen as the limit of a sequence of sharper and sharper threshold-like functions $\phi$ and $\phi^{-1}$.}
		\label{fig:threshold}
	\end{figure}

\begin{figure} [h!]
		\centering
		\includegraphics[width=2cm]{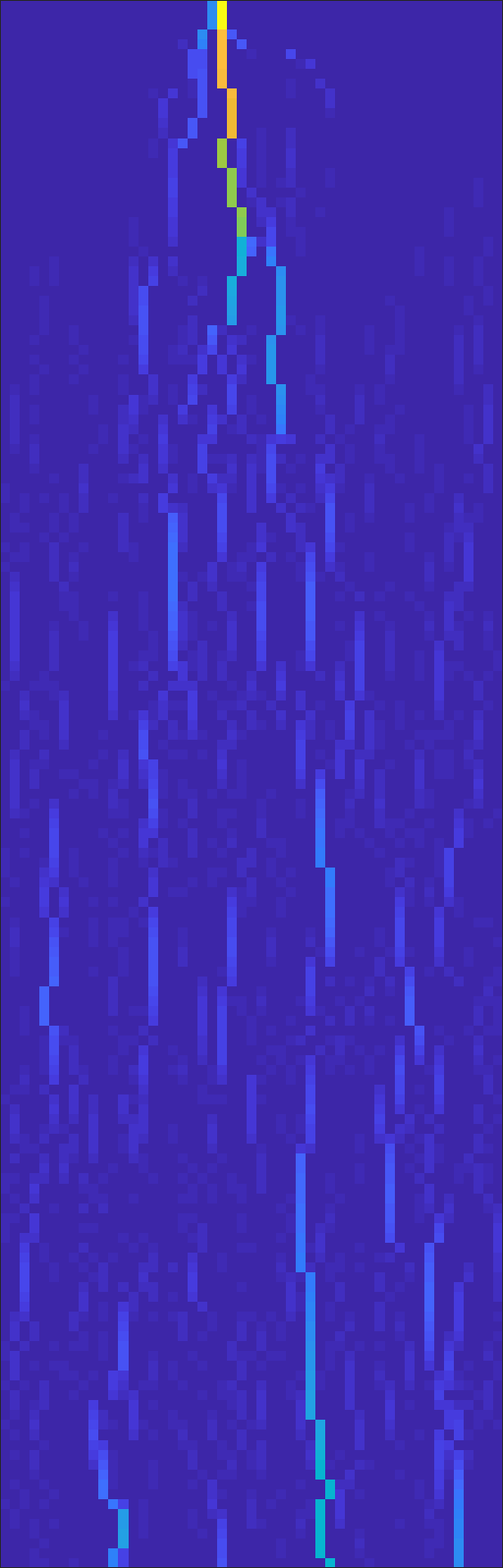}  
		\includegraphics[width=2cm]{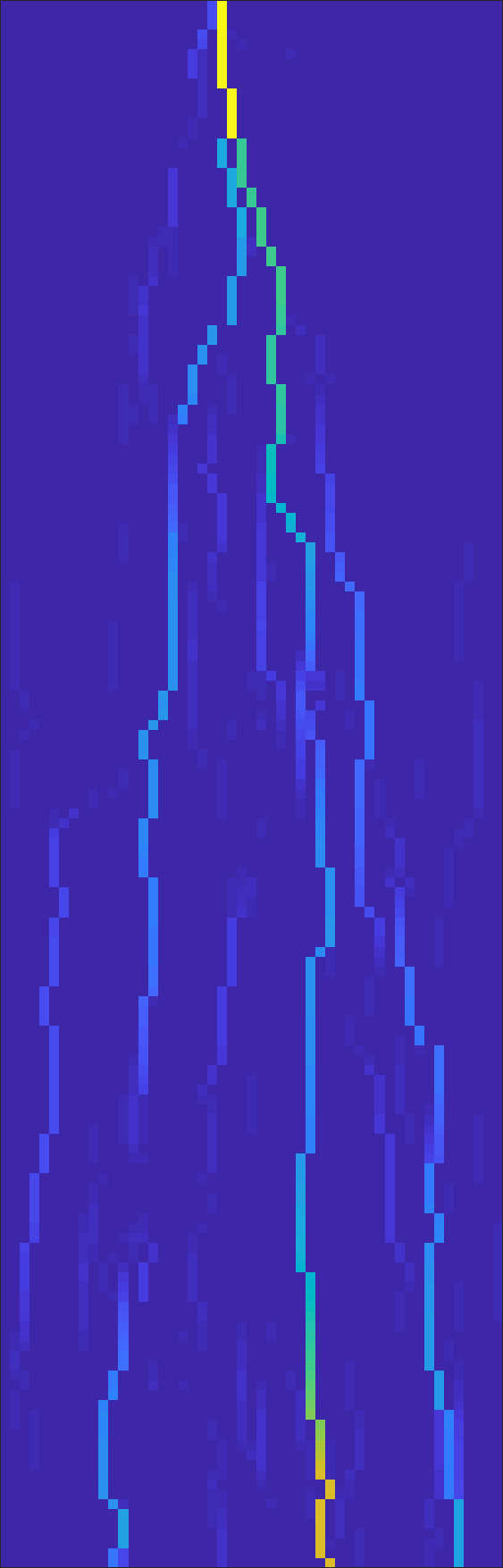}
		\includegraphics[width=2cm]{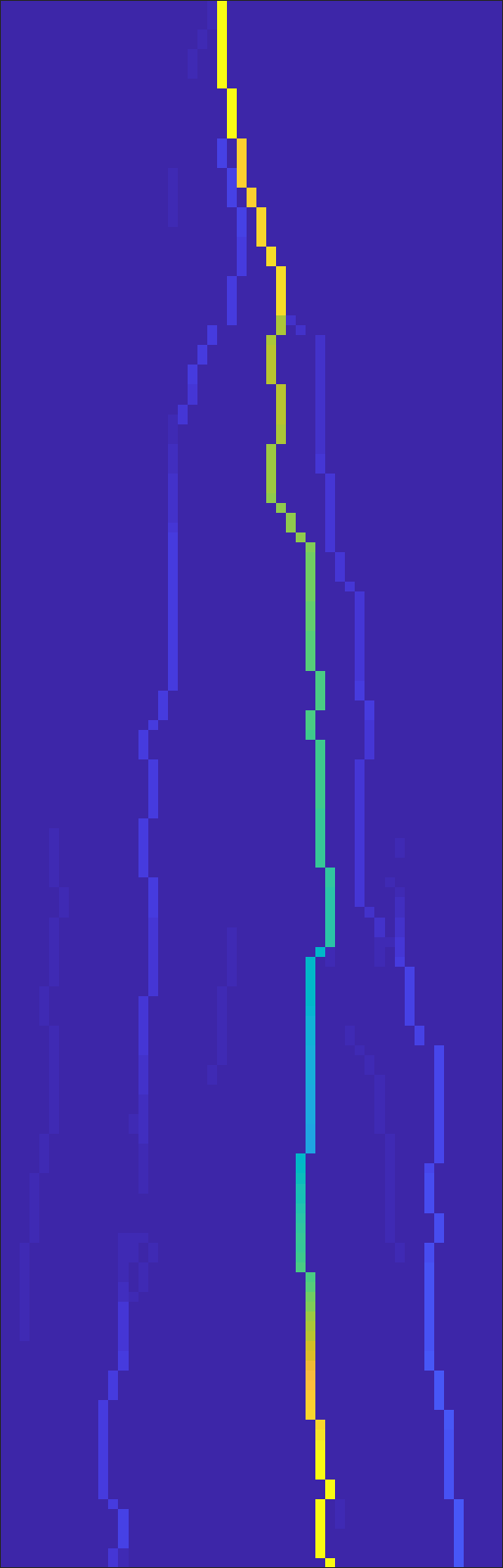}
		\includegraphics[width=2cm]{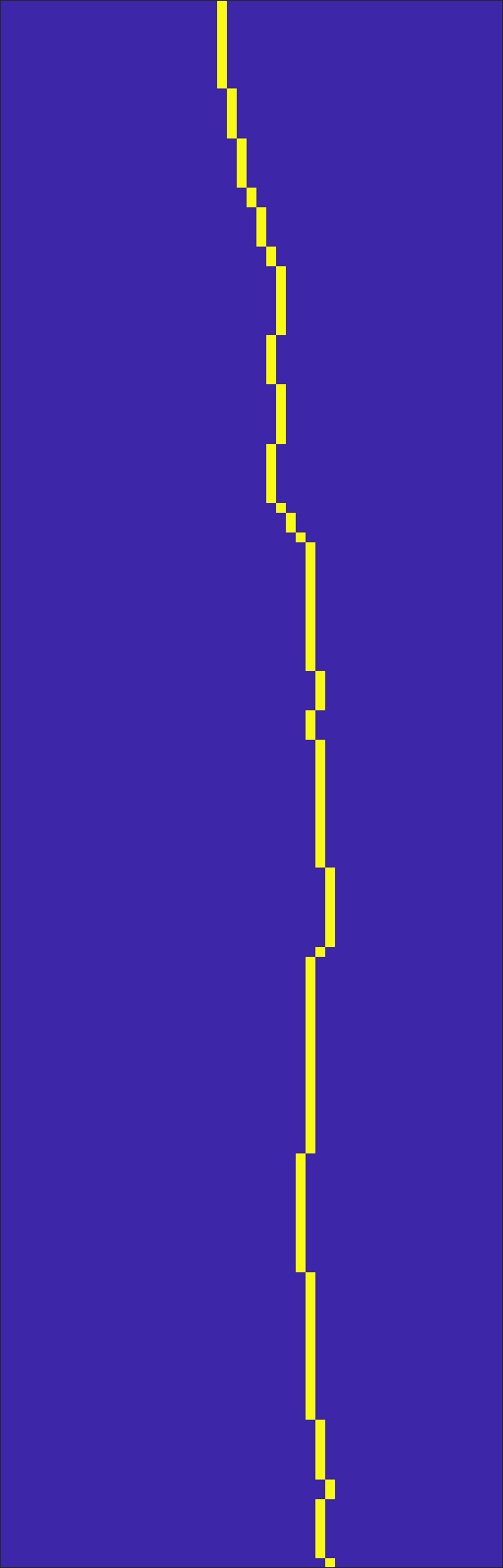}
		\caption{The lightning transient phases, from left to right: several branches are 
initially generated by the stepped leader; the stepped leader meets a connecting leader and the cloud is connected with the ground; only the main branch persists, corresponding to the optimal path (w.r.t. $J^{path}_h$); the optimal path is crossed by the long-duration current (steady state).
The color map goes from blue (no current) and light blue (low intensity current) to yellow (high intensity current).}
		\label{fig:lightning}
	\end{figure}

If the nonlinear law approaches the ideal limit characteristic, depicted in Fig. \ref{fig:threshold} (left) along with its inverse characteristics (right), then we have the following results (derived in Materials and Methods).
\begin{itemize}
\item
The steady-state current distribution minimizes the energy function
$
J^{th}= \sum_{k} V_k |u_k|.
$
\item 
Considering the family $\mathbb{P}$ of all possible paths connecting the source node $0$ to the ground node $n$,
the whole injected current $d$ flows along the path ${\mathcal P}_{h^*} \in \mathbb{P}$ that minimizes the associated total energy:
 \begin{eqnarray*}
J^{path}_h= & d \sum_{k \in {\mathcal P}_h} V_k &\rightarrow \min\\
& \mbox{s.t. }{\cal P}_h \in \mathbb{P} &
\end{eqnarray*}
where $k\in {\mathcal P}_h$ denotes that path ${\mathcal P}_h$ crosses link $k$.
\item
In the transient, the injected current starts flowing along several ``tentative" branches. When one of these branches -- corresponding to the minimum-threshold path -- first connects to the ground (in general with the aid of a connecting leader), all the others branches are depleted, as shown in Fig. \ref{fig:lightning}, obtained through simulations that faithfully reproduce this behavior.
\end{itemize}

\section*{Discussion}

Therefore, the analysis of a grid circuit with capacitors and resistors having nonlinear characteristics unravels why flow phenomena such as lightning tend to concentrate the whole current flow along a single path, despite the availability of several admissible routes: this phenomenon is due to the threshold mechanism associated with the dielectric rigidity.
In fact, for a nonlinear resistive network model, the solution of the flow equations minimizes a convex functional. In the special case of linear resistance, this functional is the dissipated energy. In the case of threshold-like nonlinear characteristics, we have proven that the minimized functional is the sum of the currents along the links, weighted by the link threshold; hence, all the current eventually flows through the global minimum path if the links connecting the grid nodes are weighed by their local dielectric rigidity.
In real situations, the dielectric rigidity can vary randomly and drastically in space, being a function of the local humidity, temperature, pressure and pollution.
This explains the seemingly random path of lightning: such a randomness is due to the gas current condition,
because the lightning  actually searches the right path.

The threshold model, including capacitive effects among nodes, faithfully describes also the transient and our simulations reproduce the behavior described in \cite{Cooray2008} and analyzed in \cite{Hager1989,Hager1998}.

Our model does not take into account inductive effects, considered by some authors in the return stroke analysis \cite[pp. 169--170]{Cooray2008}, but this does not invalidate our results because: (i) the minimum path analysis is carried out at steady state, $\dot v=0$, when the inductances are equivalent to shortcuts; (ii) the return stroke starts when the stepped leader has reached the ground, hence the path has been already ``decided''.

Our analysis reveals that lightning is one of the many phenomena in nature where a spontaneous optimization appears to take place \cite{BP2000,PD2005,Yang2014}, leading to the most efficient path choice \cite{RA19}.
In our lightning discharge model, the resulting steady-state current flow is \emph{globally} optimal, even though the current flow is \emph{locally} determined on the basis of the impedance characteristic of each single link.

\section*{Materials and Methods}

The network is modeled as a grid graph with $n+1$ nodes and $m$ links (more details on our model and our assumptions are in the SI Section \ref{app:netwmodel}). 
Node $0$ represents the node where a current $d$ is injected.
The $k$th electric component is associated with link $k$ connecting 
nodes $i$ and $j$. Its impedance is given by the parallel connection 
of a capacitance and a possibly nonlinear resistor (cf. Fig. \ref{fig:circuit} in the main paper), so that the current 
$u_k$ flowing through the component can be written as
\begin{equation}
u_k = \phi_k \Big(v_i-v_j\Big) + \frac{d}{dt}\Big[ C_k (v_i-v_j)\Big]\,,
\label{eq:admit}
\end{equation}
where $v_i$ and $v_j$ are the voltages at the terminals $i$ and $j$,
while $\phi_k(\cdot)$ is the resistor \emph{current-voltage characteristic function} and 
$C_k$ is the capacitance.

Function $\phi_k(\cdot)$ is any symmetric increasing
locally Lipschitz, or twice differentiable, function (the case in which $\phi_k(\cdot)$
is non-decreasing only is considered in the SI Section \ref{app:netwmodel}).
In the case of a linear resistor, $\phi_k=(v_i-v_j)/R_k$.
As mentioned in the main paper, we are interested in threshold-like characteristic functions
whose value is close to zero in an interval $[-V_k,V_k]$ and becomes very large if the voltage crosses the threshold value $V_k$. In Fig. \ref{fig:threshold}, the generic function $\phi$ is depicted (red curve, left) along with its inverse function $\phi^{-1}$ (red curve, right).

Following the description of Category 1 Lightning in \cite{Cooray2008,Uman2001}, we distinguish two phases.
\begin{itemize}
\item First, the stepped leader ``seeks the path to the ground'': the current is relatively low and the capacitance effect dominates, leading to a fast variation of the voltages at the nodes, in the transient evolution of the model.
\item Then, once the stepped leader has connected the cloud to the ground, dielectric breakdown is fully developed and the long-duration discharging current is triggered. We analyze this phase assuming steady-state conditions.
\end{itemize}
The initial branching phase is by far shorter than the second phase. Indeed, only by means of special very fast video equipments the first stage can be observed, while the second one can be captured by the human eye as we commonly experience. During the second phase, the current actually varies with time but its variation rate is extremely low with respect to the first phase. 

In both phases, we show that the presence of a threshold mechanism is crucial to enable the observed behavior: it explains both the transient evolution of the phenomenon and the achievement of a minimum-path steady-state configuration.

We analyze by simulations the initial transient (first phase), during which the path is ``decided" and the currents converge to a steady-state distribution. 

First of all, however, we show that in the second phase, after an initial transient, the system converges to a steady-state, i.e., terminal voltages satisfy the condition $\dot v_i = 0$.
In this state, the non-null currents define a single flow along the minimum path in terms of dielectric rigidity.

\subsection*{Steady-state analysis: the chosen path is the minimum path}
We start by showing that a threshold mechanism yields steady-state minimum-path flow.

We consider the functional
\begin{equation}\label{eq:fun}
J(u)  \doteq  \sum_{k} f_k(u_k) \doteq \sum_{k} \int_0^{u_k}~~\phi_k^{-1}(I) dI\,,
\end{equation}
where the index $k$ refers to the links (see SI Section \ref{app:netwmodel} for details).
As a first result we have the following proposition, proven in the SI Section \ref{app:proofP1}.
\begin{proposition}
\label{prop:minimum}
Given the injected current $d$, at steady state ($\dot v_i=0$)
the current distribution in the network
minimizes the functional given by \eqref{eq:fun}. If the $\phi_k$ are strictly increasing, the optimal current
distribution is unique.
\end{proposition}
Note that dimensionally $J(u_k)$ is a power, since $\phi_k^{-1}(I)$ is a tension integrated with respect to a current $I$. Consistently, in the special case of linear resistances $R_k$,  namely when $I=\phi_k(v) = v/R_k$ and $\phi_k^{-1}(I)=R_k I$,
the function in~\eqref{eq:fun} corresponds (up to the factor $1/2$) to
the  total dissipated-power distribution in \eqref{eq:dissipated}, which is the minimal
\cite[Application 1.8, Page 15]{AMO93}.

Here we are interested in the case in which the resistor characteristic is threshold-like. The ideal threshold function corresponding to a dielectric rigidity value $V_k$ is
(see Fig. \ref{fig:threshold}, blue)
\begin{equation}\label{eq:th}
\phi^{th}_k(v) =\begin{cases}
 0 & \text{if $|v| \leq V_k$}\,, \\ 
 +\infty & \text{if $v > V_k$}\,, \\ 
 -\infty & \text{if $v < -V_k$}\,.
 \end{cases}
\end{equation}
This curve is represented in blue in Fig. \ref{fig:threshold} (left).
This ideal characteristics is physically unfeasible and will not be used for our simulations. However, the corresponding optimization problem is still
well defined. Indeed,
the inverse function of $\phi^{th}_k$  (the blue curve in Fig. \ref{fig:threshold}, right) is 
\begin{equation}\label{eq:inv}
g^{th}_k(u_k) =\begin{cases}
 \text{any value in}~~ [-V_k,V_k] & \text{if $u_k=0$}\,, \\ 
 +V_k  & \text{if $u_k > 0$}\,, \\ 
 -V_k  & \text{if $u_k < 0$}\,.
 \end{cases}
\end{equation}
and for this choice the functional in \eqref{eq:fun} becomes
\begin{equation}\label{eq:func_ab}
J^{th}(u)  \doteq  \sum_{k} V_k |u_k|\,.
\end{equation}
The following proposition holds and is proven in the SI Section \ref{app:proofP2}.
\begin{proposition}
\label{prop:minimum_path}
Given the injected current $d$, the admissible (compatible with Kirchhoff's laws) 
distribution of the steady-state link currents $\bar u_k^{th}$,
$k=0,\dots,m-1$,
which minimizes functional \eqref{eq:func_ab}, corresponds to all the current
flowing from the source node to the ground along a minimum-threshold path,
namely a path ${\mathcal P}_{h^*} \in \mathbb{P}$ (where $\mathbb{P}$ is the family of paths from the source node to the ground) that minimizes the cost 
$$
J^{path} = d \sum_{k \in \mathcal{P}_h} V_k, \qquad \mathcal{P}_h \in \mathbb{P},
$$
which is the sum of all dielectric rigidities of the links along the path.
\end{proposition}
Functional \eqref{eq:func_ab} is not strictly convex, hence uniqueness is not ensured (see the SI Sections \ref{app:netwmodel} and \ref{app:proofP2} for further details).
However, the uniqueness
assumption is {\em generically satisfied}; in fact, if the 
dielectric rigidities are randomly generated, the probability of finding
two or more minimum paths with the same rigidity is zero,
hence we can assume that the minimum path is unique.

The next step is to show that the closer a characteristic function (red) is to the 
ideal threshold $\phi^{th}_k(v)$ (blue) (Fig. \ref{fig:threshold}), the closer the current
distribution is to the minimum-path distribution.
Given a sequence of characteristics $\phi^{r}_k(v)$,
$r=1,2, \dots$ which ``converge to the ideal one" and are physically feasible, so that the corresponding steady-state solutions are uniquely defined, these steady-state solutions converge to the minimum-path distribution. This property is formalized in the following theorem.

\begin{theorem}
\label{th:main}
Consider a sequence of characteristics $\phi^{r}_k$, $r=1,2,\dots$, such that, for all $r$, the corresponding steady-state current distribution $\bar u_k^{r}$ in the links
are uniquely defined. Assume that the threshold-like characteristics  converge to the ideal one
$$
\phi^{r}_k(v) \rightarrow \phi^{th}_k(v), ~~~\mbox{as}~~r \rightarrow \infty.
$$
Moreover, assume that the minimum path in terms of sum of dielectric rigidities is unique.
Then the link current distributions $\bar u_k^{r}$  converge to $\bar u_k^{th}$,
$$
\bar u_k^{r} \rightarrow \bar u_k^{th}~~~\mbox{as}~~r \rightarrow \infty\,,
$$
namely, the current distributions converge to the one with the whole current $d$ 
flowing along the minimum path.
\end{theorem}

Function $\phi^{th}_k$ is an idealized version of the gas dielectric characteristics 
in which the admittance is virtually zero for small voltage values and very large if the voltage is larger than the threshold value known as \emph{dielectric rigidity}.
In practice, true characteristics can be reasonably approximated \cite{Hager1989,Hager1998}
by a continuous curve that drastically increases after the threshold.
The meaning of the theorem is that, if these characteristics are sharp and close to $\phi^{th}_k$, then the current tends to
flow along the minimum path. The result does not rely on any specific characteristic model: only the property of the characteristics becoming ``close to the ideal"
is essential.

\begin{figure} [h!]
		\centering
		\includegraphics[width=0.3\columnwidth]{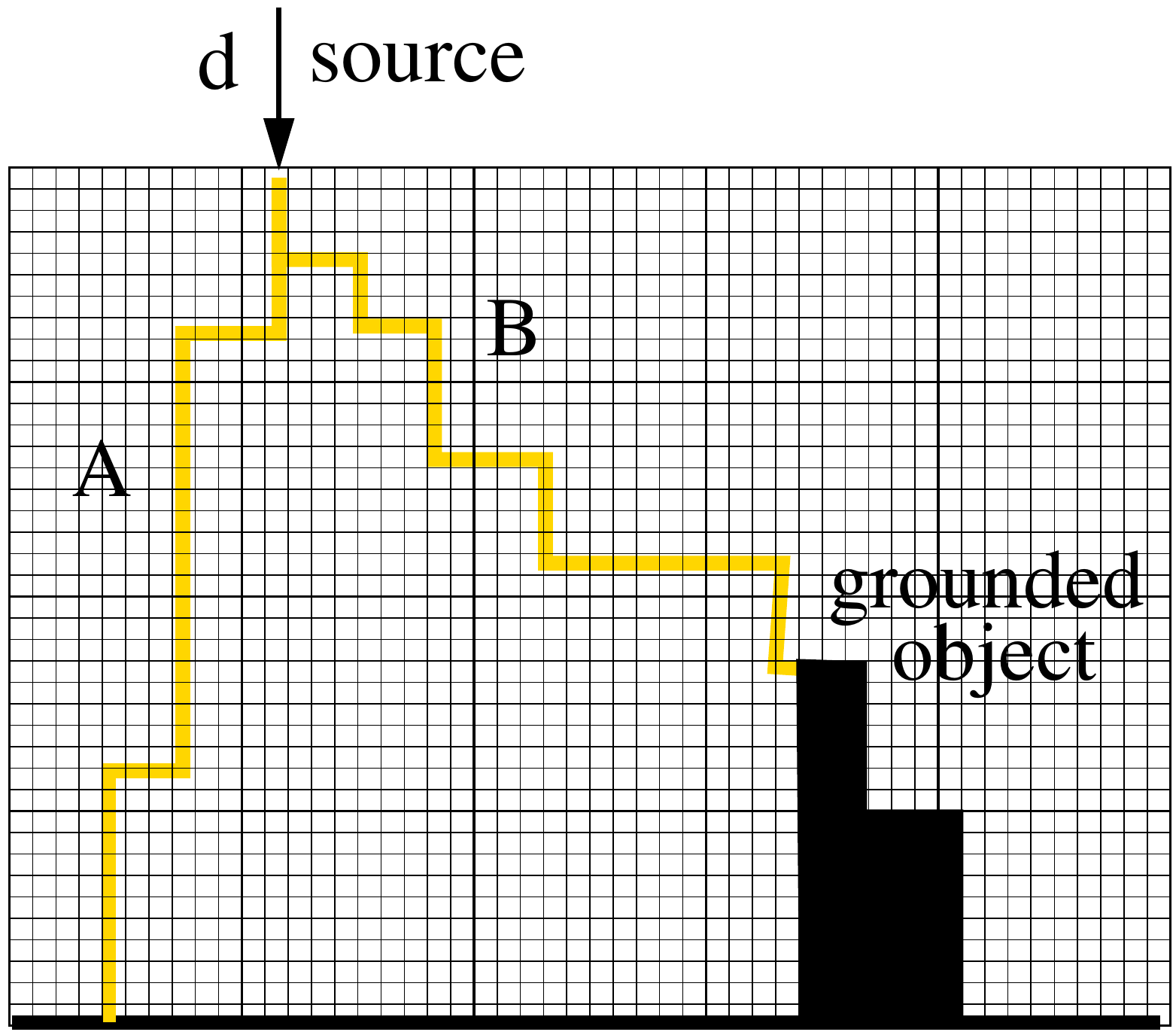} ~~~~
 		\includegraphics[width=0.25\columnwidth]{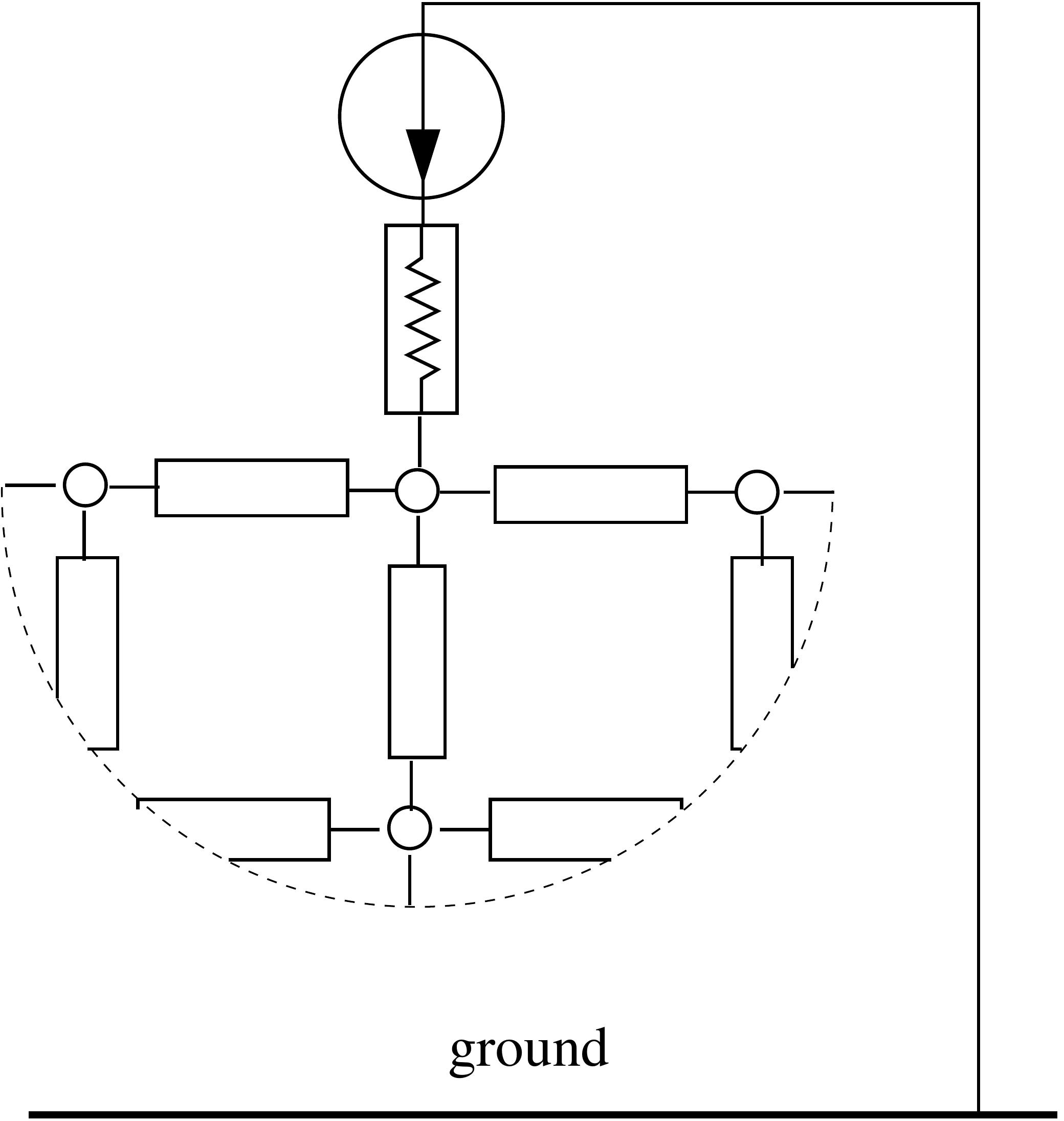}
		\caption{The model can consider, with no essential changes, different geometries.
For instance, in the case of a grounded conductive object, the cells corresponding to the object are connected by very small resistance and very large capacitance values: in the left panel, lightning would choose to reach either the ground directly (path $A$) or the grounded object (path $B$) depending on which is the minimum-threshold path, essentially considering the grounded object a zero cost portion of the path. Also voltage generators instead of current generators can be considered, as in the right panel, without any substantial change (provided that an input resistance is present).}
		\label{fig:grounded}
	\end{figure}

In the model, we consider cell-to-cell capacities,
but other capacities, such as capacities with respect to the ground, can be considered
and the analysis remains valid, since at steady state the current through the capacities is zero.
It is also fundamental to remark that the result is topology-independent: we could consider any network topology, not necessarily a grid. Also, we could consider conductive elements on the ground; in this case the lightning may find the minimum cost path
as the one that connects the source to the grounded object. Some examples are in Fig. \ref{fig:grounded}.

\subsection*{Transient analysis: seeking the minimum path}
We show here that a threshold mechanism also explains the lighting transient behavior, which can be described by the dynamic model
\begin{equation}
\label{eq:model_diff1}
BC B^\top \dot v(t) = -  B\phi \Big(B^\top v(t) \Big) + d\,.
\end{equation} 

This system asymptotically converges to the steady-state condition $\dot v(t)=0$, which leads to the condition $B \phi \Big(B^\top v(t)\Big) = d$ corresponding to the constraint of the optimization problem considered in the steady-state analysis.
A detailed stability analysis is in the SI Section \ref{app:dynamic}.

The transient analysis has been carried out via simulation, using a standard ODE solver.
In particular, to numerically demonstrate the dynamic behavior of the system, we have performed many simulations for different values of the dielectric rigidity. Videos are available to display some particularly significant cases (see the Supplementary Material Movies S1-S6 for details).\footnote{Available on-line: \url{https://users.dimi.uniud.it/~franco.blanchini/Lightsim.zip}.}

\begin{figure} [h!]
		\centering
		\includegraphics[width=0.6\columnwidth]{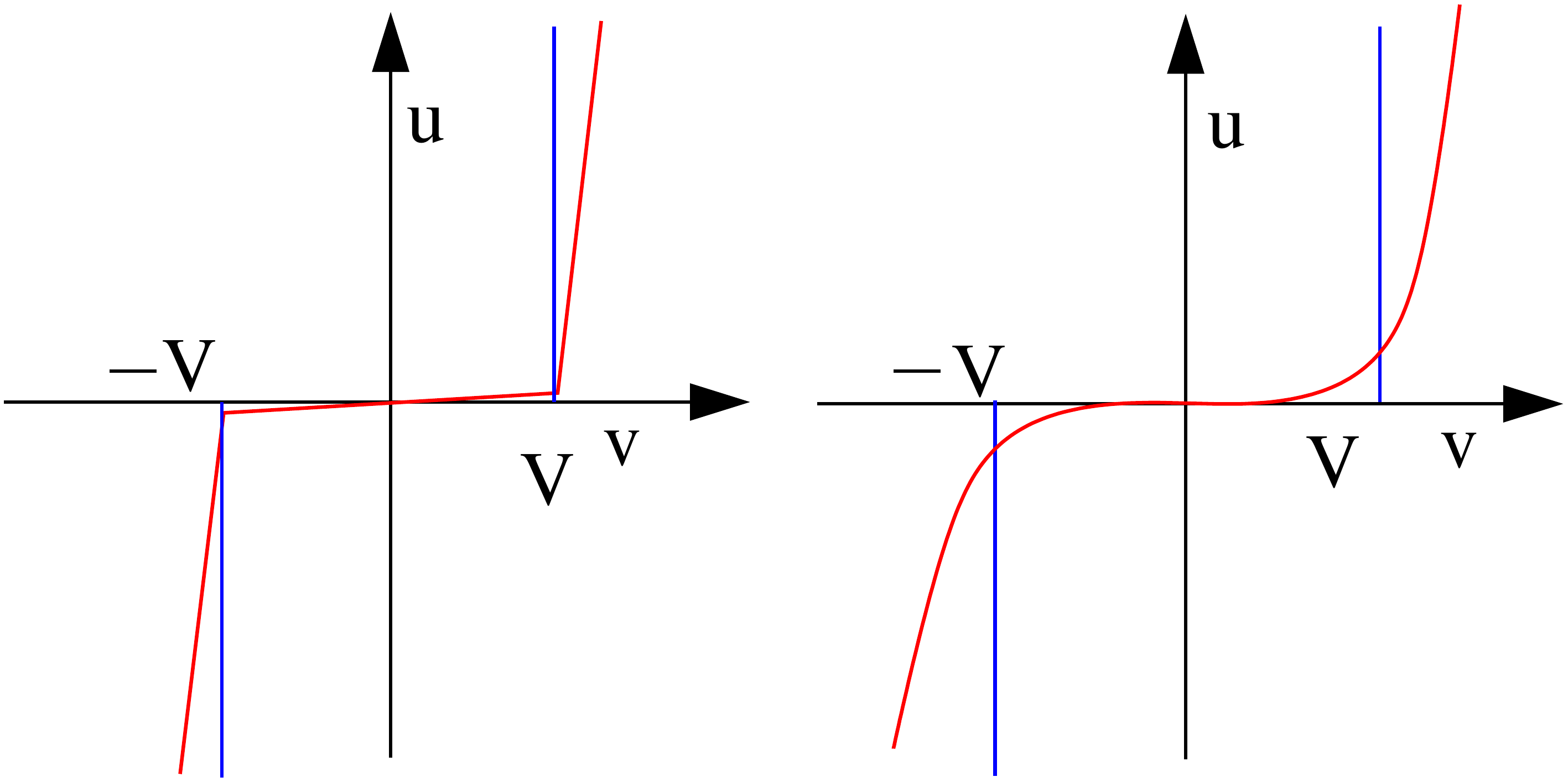} ~~~~
		\caption{The characteristics which have been used in simulation. The piecewise-linear in 
                \eqref{eq:piecewise} (left) and the monomial $(v_k/V_k)^{2r+1}$. There is no essential difference
                in the final results of the simulations. For large values of $r$, the current follows the shortest path.}
		\label{fig:characteristics_sim}
	\end{figure}

The characteristic function $\phi_k(v_k)$ can be any locally Lipschitz or continuously differentiable function that is non-decreasing and has a very high slope after the threshold.
 For simulations purposes, we have adopted the piecewise-linear threshold-like functions
\begin{equation} \label{eq:piecewise}
 \phi_k(v_k) =\begin{cases}
  \epsilon v_k & \text{if $|v_k| \leq V_k$}\,, \\ 
 r(v_k-V_k) +V_k \epsilon & \text{if $v_k > V_k$}\,, \\  
r(v_k+V_k) -V_k \epsilon& \text{if $v_k < -V_k$}\,,
\end{cases}
\end{equation}
(shown in Fig. \ref{fig:characteristics_sim}, left), all with the same plasma conductivity $r = 800$,
whereas $V_k$ has been randomly chosen for each link in the interval $V_k \in [0.5 - \delta/2, 0.5 + \delta/2]$ with uniform distribution. The variability of the dielectric rigidity is described by $\delta$, while $\epsilon$ is a very small number representing the negligible conductivity under the threshold $V_k$ (we have set $\epsilon=10^{-5}$).

Other ``sharp" characteristic functions would produce the same behavior.
We also simulated the system with the polynomial $\phi_k(v_k) =  (v_k/V_k)^{2r+1}$ (shown in Fig. \ref{fig:characteristics_sim}, right), which yields comparable results for large enough $r$, as expected. Yet, the non-Lipschitz nature of the polynomial function is numerically challenging and requires large computational times and specialized integration routines for stiff systems.

In all our numerical experiments, the capacitances have been taken all equal. We have set $C_k=1$ for all $k$, without restriction, since changing the capacitance value is equivalent to scaling time, hence the steady-state value is unaffected.
Fig. \ref{fig:lightning} reports four instants of the simulation with $\delta = 0.7$
(see Supplementary Material, Movie S5).
It can be seen that, for larger variability of the dielectric rigidity, i.e. larger $\delta$, the initial branching activity is more intense. However, the asymptotic behavior is qualitatively the same regardless of the value of $\delta$, with no exception: a single branch survives, which is numerically verified to be the minimum path in terms of total dielectric rigidity, as expected.

\subsection*{Decentralization, topology-independence and limitations}

Remarkably, as we have stressed in the main paper, the steady-state current flow in our lightning discharge model is \emph{globally} optimal, even though the current flow is \emph{locally} determined on the basis of the impedance characteristic of each single link.
In the context of distributed flow control in networks \cite{Como2017}, this kind of mechanism is called \emph{network-decentralized}~\cite{Iftar1999,ID2002,BFG15}, and localized strategies have been shown to lead to a globally optimal behavior
\cite{BFGP13,BFGMM16,Scholten2016,BCFGP19,BCFGP19b,Trip2019}.
In the considered network-decentralized control strategy, each links locally decides how much current flows through it. This approach is completely different from Dijkstra's
decentralized minimum-path algorithm, which is based on decentralized
dynamic programming techniques \cite{Luenberger,ZhangLi2017} and in which the routing decision is made at the nodes: each node locally decides to which of the outgoing links an incoming unit of flow must be redirected.
In our setup, a ``link decision'' is made: each admittance can be interpreted as an agent
that locally decides how much current is allowed to flow.

Our results are independent of the topology of the network, which does not necessarily need to be a grid graph with square cells: other topologies would lead to a minimum path solution.

However, we stress that our model is far from capturing all the complex aspects of lightning.
Its validity is limited to the beginning of the phenomenon, until the return stroke is triggered,
because in this initial stage the path is chosen. 
After the stepped leader has reached the ground, the current
follows  the ``chosen minimum path" until the end, as it is experimentally
well documented (and confirmed by our simulations), because this path becomes a ionized  
channel with low resistance. 
So our model and simulations are not expected to be a faithful quantitative reproduction of the whole lightning process (including discharge endurance, power
dissipation and so on), but their significance is limited to the first part. Yet, 
the qualitative behavior, with {\em the discharge following the minimum path},
has been always confirmed with no exception.
Moreover, our model does not consider other aspects such as the ground currents, which are not relevant to the path choice. 
We have considered the so-called Category 1 Lightning, cloud-to-ground,
which is the most common type of lightning; however, the model can be adapted to any type of lightning of gas discharge.

\newpage

\section{Supplementary Information}

\subsection{Network model and assumptions}\label{app:netwmodel}

The electrical network is associated with a graph $\mathcal{G}=(\mathcal{N},\mathcal{L})$,
where $\mathcal{L} =\{0, \ldots, m-1\}$ is the set of the $m$ links, each modeling an electric component, and $\mathcal{N}=\{0,\ldots, n\}$ is the set of the $n+1$ nodes, each modeling a terminal where some components join. 
In particular, node $n$ corresponds to the zero-potential ground, while at node $0$ a (current or voltage) generator is applied, with its other terminal grounded, inducing an input current $d$ that enters the network.  

Consider the $k$-th electric component of the network, associated with the graph link $k \in \mathcal{L}$, $k = (i,j)$, which connects node $i$ to node $j$: its impedance is given by the parallel connection of a capacitance and a possibly nonlinear resistor, so that the current $u_k$ flowing through the component can be written as
$u_k = \phi_k \Big(v_i-v_j\Big) + \frac{d}{dt}\Big[ C_k (v_i-v_j)\Big]$ (corresponding to the admittance equation~\eqref{eq:admit}),
where $v_i$ and $v_j$ are the voltages at the terminals $i$ and $j$,
while $\phi_k(\cdot)$ is the resistor \emph{current-voltage characteristic function} and 
 $C_k$ is the capacitance.

We consider the \emph{generalized node-link incidence matrix} of the graph $\mathcal{G}$, which is the matrix $B \in \{-1,0,1\}^{n \times m}$ obtained by assigning an arbitrary direction to each link $k = (i,j)$ of $\mathcal{G}$ and setting a $1$ entry in position $i$ (source node) and a $-1$ entry in position $j$ (destination node), and zero elsewhere, in the corresponding $k$-th column of $B$, and then removing the row corresponding to node~$n$. In particular, the $m$ columns of $B$ are associated with the links representing the electric components and its $n$ rows are associated with the nodes representing terminals. Links coming from the external environment (associated for instance with an injected current) have a single nonzero entry, equal to $-1$, corresponding to their destination node and links going to the external environment have a single nonzero entry, equal to $1$, corresponding to their source node: in our model, the connections to the external environment are the (zero-voltage) ground and the source of supplied power.
The following assumption holds, because we have considered the ground zero-potential node as an external node.

\begin{assumption}
The network graph $\mathcal{G}$ is connected internally and connected to the external environment.
As a consequence, matrix $B$ has full row rank.
\end{assumption}

Links are associated with the currents $u_0,\ldots,u_{m-1}$ flowing through the individual electrical components, which we group in the vector $u\in \mathbb{R}^m$; nodes are associated with terminal voltages $v_0,\ldots,v_{n-1}$, grouped in the vector $v\in \mathbb{R}^n$.

Then, the admittance equation~\eqref{eq:admit} can be written as
\begin{equation}
u_k = \phi_k \Big(B_k^\top v \Big) + C_k B_k^\top \dot v\,,
\label{eq:admit1}
\end{equation}
while the current balance at the node $i$ is
\begin{equation}
B^i u-d_i=0\,,
\label{eq:balance}
\end{equation}
where $B^i$ is the $i$-th row of $B$ and $d_i$ is the $i$-th element of the vector of externally supplied current.

Merging equations \eqref{eq:admit1} and \eqref{eq:balance} yields the dynamics of the overall circuit $\mathcal{G}$, in terms of voltages and currents, which is described by the discretized space model \cite{Hager1989,Hager1998} with equations
\begin{equation}\label{eq:model_diff}
 0 = Bu(t) - \bar d, \quad u(t) = \phi \Big(B^\top v(t) \Big) + C B^\top \dot v(t)\,,
\end{equation}
where $C=\mbox{diag}\{C_0,C_1,\dots, C_{m-1}\}$ is a diagonal matrix whose diagonal elements are the capacities $C_k$, 
$u=[u_0,u_1,\dots, u_{m-1}]^\top$ is the vector whose components are the currents along the links of $\mathcal{G}$, 
$\bar d=[d,0,\dots,0]^\top \in \mathbb{R}^n$ is the input current vector,
$v=[v_0,v_1,\dots, v_{n-1}]^\top$ is the vector whose components are the voltages at the nodes of $\mathcal{G}$, $\dot v$ is the time derivative of vector $v$ and $\phi(\cdot) =[\phi_0(\cdot),\phi_1(\cdot),\dots,\phi_{m-1}(\cdot)]^\top$ is the vector of the characteristic functions.

The following general  assumption is considered.

\begin{assumption}\label{ass:phi} 
	Each characteristic function
	$\phi_k: \mathbb{R}\rightarrow \mathbb{R}$ is a possibly nonlinear odd 	monotonically increasing function
and locally Lipschitz.
\end{assumption}

This assumption can be weakened by requiring $\phi_k$ to be \emph{monotonically non-decreasing} only.

The shape of a generic characteristic function satisfying Assumption~\ref{ass:phi} is shown in red in Fig.~\ref{fig:threshold}. The assumption implies that, for each link $k \in \mathcal{L}$,  $\phi_k$ is invertible.
Then, for each $k$, the function
\begin{equation}\label{eq:f}
f_k : y \mapsto \int_0^y g_k(s)ds\,,
\end{equation}
where $g_k \doteq \phi_k^{-1}$ is the (monotonically increasing) inverse function of $\phi_k$, is well defined in $(-\infty,+\infty)$.

Functions $f_k$ are continuously differentiable. In addition, they are strictly convex, since their derivative $g_k$ is an increasing function almost everywhere ($f_k'' = g_k' >0$ almost everywhere; in fact, $g_k^\prime$ may be not defined in some isolated points, e.g. in $u=0$ for $g_k(u) = u^{1/3}$).

If we assume that $\phi_k$ is non-decreasing only, then we have convexity but not strict convexity of $f_k$.

\subsection{The dynamic model}\label{app:dynamic}
We report here the stability analysis of the complete model \eqref{eq:model_diff},
which we can rewrite in the equivalent form
\begin{equation*}
	\dot v(t) =  - [BCB^\top]^{-1}\left[ B\phi \Big(B^\top v(t) \Big) - \bar d\right ]\,.
\end{equation*}
The stability of this type of systems has been studied in the literature
\cite{VDS_SCL,BCFGP19}. Consider the steady-state vector~$\bar v$, such that 
\begin{equation*}
0 =  [BCB^\top]^{-1} \left [B \phi \Big(B^\top \bar v \Big) - \bar d\right ] \,,
\end{equation*}
and denote by $x$ the shifted variable defined as $x(t)=v(t)-\bar v$, whose time variation is
\begin{equation*}
  \dot x(t) =  [BCB^\top]^{-1}B\left [\phi \Big(B^\top (x(t)+\bar v) \Big) -\phi \Big(B^\top \bar v \Big)\right ] \,.
\end{equation*}
Since $\phi$ is a vector of strictly increasing functions, we can write
\begin{equation*}
\phi \Big(B^\top (x+\bar v) \Big) -\phi \Big(B^\top \bar v \Big)) = \Delta(v(x)) B^\top x(t)\,,
\end{equation*}
where $\Delta(v)$ is a diagonal matrix of strictly positive continuous functions \cite{BFGMM16,BCFGP19,BCFGP19b}. Hence
\begin{equation*}
\dot x =  -[BCB^\top]^{-1}B\Delta(v(x)) B^\top  x\,.
\end{equation*}
Consider the positive definite Lyapunov function candidate $V(x)=\frac{1}{2} x^\top BCB^\top x$, which is the energy stored in the capacitors. Its derivative is negative definite:
\begin{equation*}
\dot V(x) = x^\top BCB^\top \dot x = - x^\top B\Delta(v) B^\top  x <0
\end{equation*}
as $x \neq 0$. This ensures asymptotic stability of the steady-state solution.

\subsection{Proof of Proposition \ref{prop:minimum}}\label{app:proofP1}

We have to prove that the steady-state current distribution in the network induced by a constant current injection $\bar d$, achieved when $\dot v =0$, namely when
\begin{equation}
B \phi \Big(B^\top v \Big)- \bar d=0\,,
\label{eq:solving_st}
\end{equation}
is indeed the current that solves the optimization problem
\begin{eqnarray} \label{eq:funct}
J(u_0,\ldots,u_{m-1}) &\doteq& \sum_{k=0}^{m-1} f_k(u_k) \rightarrow \min \\
\mbox{s.t.}&&  Bu = \bar d \,, \label{eq:const} 
\end{eqnarray}
where \eqref{eq:const} is the flow constraint imposed by Kirchhoff's current law
\cite{AMO93}.

In view of the assumptions on $\phi_k$, the function $f_k$ defined in \eqref{eq:f} is continuously differentiable with
increasing derivative, hence strictly convex.
Therefore, the optimization problem~\eqref{eq:funct}-\eqref{eq:const} is strictly 
convex and has a unique solution, achieved by applying the first order Karush-Kuhn-Tucker conditions to the Lagrangian function
\begin{equation*}
\sum_{k=0}^{m-1} f_k(u_k) + \lambda^\top [Bu - \bar d]\,,
\end{equation*}
where $\lambda \in \mathbb{R}^n$ is the vector of Lagrangian multipliers.
The derivative with respect to $u$ must be zero, hence we get
\begin{equation}
\label{eq:lagran}
\nabla f(u) + \lambda^\top B = 0\,.
\end{equation}
Now, the first derivative of the elements of $f$ is $f'_k=g_k$, which is invertible with inverse $\phi_k$. As a consequence, $u_k = \phi_k \Big(B_k^\top \lambda \Big)$, hence
 \begin{equation}
\label{eq:lagran1}
u= \phi \Big(B^\top \lambda \Big)\,. 
\end{equation}
The solution of the optimization problem is therefore the unique solution of the system \eqref{eq:solving_st}, $B \phi \Big(B^\top \lambda \Big) - \bar d =0$. Interestingly, the Lagrange multiplier vector is the steady-state voltage, $\lambda=v(\infty)$.

In the case of non-decreasing functions $\phi_k$, we still have convexity but not strict convexity: the result holds, but the minimizing distribution may be non-unique.

Finally note that, in the special case of linear resistances $R_k$, namely when $u_k = (v_i-v_j)/R_k$,
the solution of the optimization problem provides (half) the minimum-dissipated-energy distribution, hence 
\begin{equation*}
J^{\star}= \frac{E_{tot}}{2}=\frac{1}{2}\sum_{k=0}^{m-1} R_k u_k^2\,;
\end{equation*}
see for instance \cite[Application 1.8, Page 15]{AMO93}.
However, in the general case, each $f_k$ can be different from the local dissipated energy, which is $\frac{1}{2}E_{k}=\frac{1}{2}u_k\phi_k^{-1}(u_k)$.
Hence, in the nonlinear case, the minimized functional is not the dissipated energy as in the case of linear resistances.

\subsection{Proof of Proposition \ref{prop:minimum_path}}\label{app:proofP2}

We have to show that the limit optimization problem
\begin{eqnarray} \label{eq:funct1}
J^{th}(u) &\doteq& \sum_{k=0}^{m-1}f^{th}_k(u_k) = \sum_{k=0}^{m-1} V_k|u_k| \rightarrow \min \\
\mbox{s.t.}&& Bu = \bar d\,, \label{eq:funct2}
\end{eqnarray}
associated with the limit characteristic $\phi^{th}_k$ (the limit of $\phi^{r}_k$ as $r \rightarrow \infty$),
admits as its optimal solution the current distribution with all current $d$ channeled along the shortest path.

Assume that the injected current is positive: $d>0$ (the case $d<0$ is identical).
Let $u^*$ denote the current distribution solving the optimization problem \eqref{eq:funct1}--\eqref{eq:funct2}.
To keep the proof simple we assume that all links in the network have been oriented in such a way that
$u_{k} \geq 0$. This is not a restriction since link orientation is arbitrary.

Now consider the modified problem
\begin{eqnarray} \label {eq:funct3}
 & & \sum_{k=0}^{m-1} V_k u_k  \rightarrow \min \\
\mbox{s.t.}&& Bu = \bar d\,, \label{eq:funct4} \\
              && u \geq 0 \,, \label{eq:funct5} 
\end{eqnarray}
where the absolute value has been removed and a positivity constraint has been added.

The solution $u^*$ of the previous problem is a feasible solution of the new problem,
because its elements are nonnegative by construction.
It is also optimal for the new problem. Indeed, if another solution $\tilde u$
were found with a lower cost, this would be a feasible solution also for the original problem
\eqref{eq:funct1}--\eqref{eq:funct2} and would have a cost smaller than that of $u^*$.

Now note that, to solve \eqref{eq:funct3}--\eqref{eq:funct5},
we can just take $d=1$ and then scale the solution (by the true value $d >0$).
The proof is concluded by noticing that \eqref{eq:funct3}--\eqref{eq:funct5}
with $d=1$ gives the minimum cost path \cite{AMO93}, with optimal cost
$d \sum_{k=0}^{m-1} V_k$
and the whole flow through the minimum cost path.
An energetic interpretation is that the solution minimizes the overall power,
measured as the product between the current flowing in a link and its dielectric rigidity.

\subsection{Proof of Theorem \ref{th:main}}\label{app:proofT1}

We have to prove that the steady-state solutions $u^{r*}$ associated with the characteristic 
functions $\phi^{r}_k$, which have been shown to be the minimizers of
\begin{eqnarray} \label {eq:new1}
J(u_0,\ldots,u_{m-1}) &\doteq& \sum_{k=0}^{m-1} f^r_k(u_k) \rightarrow \min \\
\mbox{s.t.}&&  Bu = \bar d \,, \label{eq:n} 
\end{eqnarray}
converge to the solution of \eqref{eq:funct1}--\eqref{eq:funct2}
if this is unique  (equivalently, the minimum path is unique).

Denote by $J^r$ and $J^{th}$ the cost functionals of the considered optimization problems,
\begin{equation*}
J^r(u) = \sum_{k=0}^{m-1} f^r_k(u),~~J^{th}(u)=\sum_{k=0}^{m-1}  V_k  \left  | u_k \right |\,.
\end{equation*}
Since $g_k^r$ is strictly increasing, $f_k^r$ is strictly convex and, in turn, also $J^r$ is strictly convex. Hence, the minimizer vector $u^{r*}$ of \eqref{eq:funct}-\eqref{eq:const} is unique.

Let $J^{th*}=J^{th}(u^*)$ be the optimal cost of the limit problem.
For any $y$, $f^r_k(y)\underset{r \rightarrow \infty}{\longrightarrow} V_k  \left  | y \right |$.
Then, $J^r(u) \underset{r \rightarrow \infty}{\longrightarrow} J^{th}(u)$ and, in particular, $J^r(u^*) \underset{r \rightarrow \infty}{\longrightarrow} J^{th}(u^*)= J^{th*}$.
This means that  the sequence of optimal costs $\{J^r(u^{r*})\}_{r\in\mathbb{N}}$ is upper bounded by a sequence $\{J^r(u^*)\}_{r\in\mathbb{N}}$  that converges to $J^{th*}$ as  
 \begin{align}
J^r(u^{r*}) \leq J^r(u^*) \underset{r \rightarrow \infty}{\longrightarrow} J^{th*}\,. \label{eq:bou}
 \end{align}
Functionals $J^r$, for all $r$, as well as $J^{th}$, are radially unbounded because they are the sum of non-negative radially unbounded functions $f^r_k$.
Then, their optimal solutions $u^*$, respectively $u^{r*}$ are finite.
In view of \eqref{eq:bou} there exists $\bar J >0$ for which  these solutions $u^{r*}$
are inside the compact set
\begin{equation*}
 \mathcal{S}^{r*} = \{u\in\mathbb{R}^m\,:\,J^{r*}(u)\le \bar{J}\}\,.
\end{equation*}
By construction $J^r(u)$ converges to $J^{th}(u) = \sum V_k |u_k|$, which is radially unbounded. Then, the sequence of  sets $\mathcal{S}^{r*}$ is uniformly bounded in a compact set $\mathcal{S}$, hence all optimal solutions are uniformly bounded: $\{u^{r*}\}_{r\in\mathbb{N}}  \in  \mathcal{S}$.

We prove the convergence $u^{r*}  \underset{r \rightarrow \infty}{\longrightarrow} u^*$ by contradiction.
We assume that  $u^{r*} \underset{r \rightarrow \infty}{\not\longrightarrow} u^*$.
Negating convergence to $u^*$ implies that there exist an open neighborhood, $\mathcal{U} \subset \mathcal{S}$, of $u^*$ and a sub-sequence of $u^{r*}$ that is in the complement of $\mathcal{U}$ in $\mathcal{S}$, namely in the compact set $\mathcal{S} \setminus \mathcal{U}$.
In turn, this sub-sequence confined in the compact set admits a sub-sub-sequence that converges to some point $u^\circ \in \mathcal{S} \setminus \mathcal{U}$.
Hence there exist a sub-sequence  
$\{u^{\tilde{r}*}\}_{\tilde{r}\in\mathcal{N}}$ of the original  $\{u^{r*}\}_{r\in\mathbb{N}}$ that converges to some vector $u^\circ \neq u^*$, 
being $\mathcal{N}=\{N_1,N_2,\dots \}$ an infinite ordered set of increasing integers.
All vectors $u^{\tilde{r}*}$ satisfy the constraint $Bu^{\tilde{r}*}=\bar d$ as they are solutions to problem~\eqref{eq:funct}--\eqref{eq:const}.

Hence, also the limit vector $u^\circ$ does satisfy
$B u^\circ=\bar d$. 
Then, the proof can be concluded by showing that 
 \begin{align}
J^{th}(u^\circ) \leq J^{th*} \,, \label{eq:contrad}
 \end{align}
which is a contradiction, because it would imply that either $J^{th*}$ is not the optimal as assumed, or (if equality holds) that the optimization problem~\eqref{eq:funct1} has two minimum points, $u^*$ and $u^\circ$, against the uniqueness assumption.

To prove \eqref{eq:contrad}, the first step is to note that there exists a finite value $g$ such that, for all $k$ and for all sufficiently large $r$, $g^r_k(u_k) \leq g$ for all $u_k$ such that $u \in \mathcal{S}$, since $g^r_k(u_k) \underset{r \rightarrow \infty}\longrightarrow V_k$ and $\mathcal{S}$ is a compact set, and hence is bounded.
As a consequence, for a sufficiently large $r$, $J^r$ has a uniformly bounded gradient since, for all $y$,
\begin{equation*}
\left |\partial J^r /\partial u_k  \right |(y)= g^r_k(y)\leq \max_{u_k:u \in \mathcal{S}}\{g^r_k(u_k)\} \leq g\,.
\end{equation*}
Hence, as $J^r$ is convex for all $r$, there exists a constant $L$ such that, for sufficiently large $r$ and for all $p,q\in\cal{S}$,
\begin{equation*}
\left | J^r(p)- J^r(q)\right | \leq L \|p-q\|\,.
\end{equation*}
Then, for  $\tilde{r} \in \mathcal{N}$ sufficiently large, the following inequalities hold: 
\begin{align*}
\hspace{-5mm} J^{th}(u^o) &= J^{th}(u^o) - J^{\tilde{r}}(u^o) + J^{\tilde{r}}(u^o) - J^{\tilde{r}}(u^{\tilde{r}*}) +J^{\tilde{r}}(u^{\tilde{r}*})\\
 &\leq \underbrace{|J^{th}(u^o) - J^{\tilde{r}}(u^o)|}_{\rightarrow 0} +
\underbrace{|J^{\tilde{r}}(u^o) - J^{\tilde{r}}(u^{\tilde{r}*})|}_{\leq L \| u^o -u^{\tilde{r}*}\|\rightarrow 0} + J^{\tilde{r}}(u^{\tilde{r}*}) \\ 
 &\leq J^{\tilde{r}}(u^*) \underset{\tilde{r} \rightarrow \infty}\longrightarrow J^{th*}\,,
\end{align*}
where the last inequality and the limit come from \eqref{eq:bou}
and the fact that $\{u^{\tilde{r}*}\}_{\tilde{r}\in\mathcal{N}}$  is a sub-sequence of $\{u^{r*}\}_{r\in\mathbb{N}}$ with limit $u^o$. Then we have shown the contradiction \eqref{eq:contrad}, which concludes the proof.


\end{document}